\documentclass[10pt,twocolumn,letterpaper]{article}

 \usepackage{cvpr}              

\usepackage{comment}

\usepackage{soul}
\setuldepth{foobar}

\usepackage{multirow}
\usepackage{arydshln}
\usepackage{listings}
\usepackage[most]{tcolorbox}
\usepackage[T1]{fontenc}
\usepackage{textcomp}

\lstdefinestyle{ieeeprompt}{
  basicstyle=\ttfamily\scriptsize,
  columns=fullflexible,
  breaklines=true,
  breakatwhitespace=true,
  keepspaces=true,
  showstringspaces=false,
  upquote=true,
  aboveskip=0.5\baselineskip,
  belowskip=0.5\baselineskip
}

\newtcblisting{promptbox}[2][]{%
  enhanced,
  arc=1mm, 
  breakable,   
  boxrule=0.4pt,
  colback=white,
  colframe=black!95,
  titlerule=0pt,    
  toptitle=0mm,        
  bottomtitle=0mm,     
  top=0.6mm,           
  left=1.2mm,right=1.2mm,top=0.8mm,bottom=0.8mm,
  listing only,
  listing options={style=ieeeprompt},
  title=\textbf{#2},
  fonttitle=\footnotesize,
  #1
}

\usepackage{gensymb}

\definecolor{beige}{RGB}{255, 250, 240}  

\usepackage{framed}

\definecolor{cvprblue}{rgb}{0.21,0.49,0.74}
\usepackage[pagebackref,breaklinks,colorlinks,allcolors=cvprblue]{hyperref}

\def\paperID{*****} 
\def\confName{CVPR}
\def\confYear{2026}

\title{Semantic Constraint Synthesis for Adaptive Trajectory Optimization via Large Language Models}

\author{Eleanor Brosius, Yuji Takubo, Daniele Gammelli, Simone D'Amico, Marco Pavone\\
Stanford University \\
496 Lomita Mall, Stanford, CA, 94305\\
{\tt\small \{elbro, ytakubo, gammelli, damicos, pavone \}@stanford.edu}}
\begin{document}
\maketitle

%
%
%
\begin{abstract}
Trajectory optimization is a critical component for enabling safe and reliable autonomous operations in space exploration. As space missions increase in frequency, complexity, and scope, there is a growing need to rapidly formulate mathematically sound trajectory optimization problems that accurately reflect mission objectives and operational constraints. However, translating mission intent into tractable analytical formulations for trajectory optimization requires substantial domain expertise. This paper presents a framework that leverages large language models (LLMs) to translate natural language descriptions of mission requirements and constraints into executable trajectory optimization code and corresponding mathematical formulations. Experiments in spacecraft rendezvous scenarios demonstrate a high success rate in reconditioning a convex trajectory optimization problem from semantic mission requirements. Ultimately, this work highlights the potential of LLMs to bridge high-level intent and formal optimization models, enabling more flexible and efficient trajectory design of spacecraft.

\end{abstract}

\section{Introduction}

\begin{figure*}[t!]
    \centering
    \includegraphics[width=1.0\linewidth]{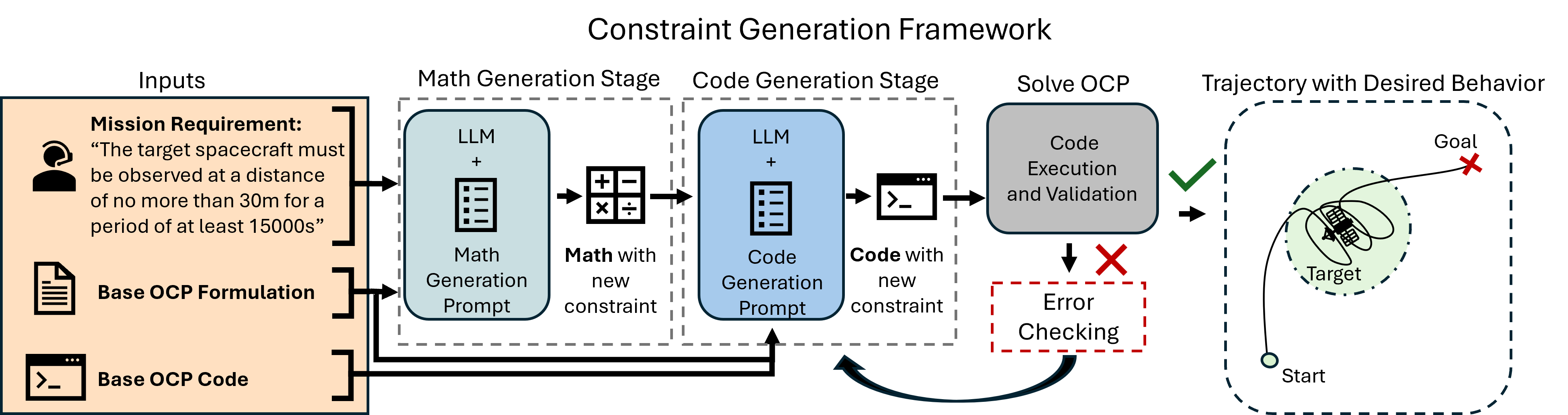}
    \caption{\textbf{Automated constraint synthesis from natural-language mission requirements.} A structured input consisting of (i) the baseline \LaTeX{} formulation of the OCP, (ii) the corresponding baseline code, and (iii) a text-based mission requirement is provided as input to the proposed framework to initialize a structured generation process. First, a pretrained LLM translates the requirement into a formal mathematical expression that augments the original problem formulation. Next, a subsequent query generates the corresponding executable code snippet. Finally, the updated code is solved by the numerical OCP solver, producing an adaptive and safety-compliant spacecraft trajectory.}
    \label{fig:main}
\end{figure*}

The formulation of optimal control problems (OCPs) plays a critical role in spacecraft trajectory generation, enabling both high-level mission planning and low-level maneuver execution. 
Although numerous case-specific solution methods have been developed by exploiting particular vehicle dynamics \cite{dsouza1997optimal, chernick_phd_2021}, there is an increasing demand for general-purpose solvers capable of handling constrained trajectory optimization in complex scenarios such as rocket landing and space rendezvous.
Convex optimization, in particular, provides a powerful framework due to its strong theoretical foundations \cite{boyd2004convex, malyuta_scp_2022} and the availability of efficient embedded solvers \cite{domahidi2013ecos}, making it well-suited for real-time and safety-critical applications.
However, its practical implementation and iterative design are often constrained by the need for substantial domain expertise and labor-intensive problem reformulation. Moreover, in operational settings where environmental conditions or mission objectives may change, reliance on a fixed, parameterized OCP significantly limits system adaptability.

In parallel, recent advances in large language models (LLMs), capable of scenario interpretation, mathematical reasoning, and code generation, have enabled the translation of natural-language specifications into formal, solver-verifiable constraints and symbolic representations \cite{ramamonjison2023nl4opt,peng2026nc2c}. Despite the existing theoretical and algorithmic support for embedding complex constraints within a convex-based OCP framework\cite{szmuk2020successive, malyuta2023fast, malyuta_scp_2022}, current spacecraft guidance methodologies have not yet capitalized on these advances. Leveraging these algorithms would provide a principled and systematically extensible foundation for language-driven problem formulation.


This paper proposes an LLM-based semantic interface that compiles natural-language mission requirements into solver-ready trajectory constraints for OCPs. Given a baseline problem specification and a newly documented mission requirement, the proposed pipeline generates both a formal mathematical formulation of the corresponding constraints and executable Python code that augments the underlying OCP.
We situate this contribution within the broader context of aerospace vehicle guidance, which traditionally demands substantial expertise in vehicle dynamics and optimal control. Incorporating language-driven problem formulation offers practical advantages by enabling more intuitive specification of guidance objectives while streamlining iterative design and refinement.

The contributions of this work are threefold:
\begin{itemize}
    \item We introduce a semantic interface that translates natural-language mission requirements into formal mathematical constraints and executable code compatible with numerical optimal control solvers using LLMs
    \item We conduct initial ablation studies to analyze the impact of different LLM prompting and generation strategies, providing insight into the contributions of individual components within the semantic compilation pipeline.
    \item We validate the proposed framework in a preliminary spacecraft rendezvous case study, demonstrating its ability to reformulate and augment a convex optimal control problem based on semantic constraint descriptions.
\end{itemize}

\section{Related Work}

\noindent \textbf{Semantics-based Adaptive Trajectory Generation: } 
The majority of learning-based spacecraft guidance methods target efficient solutions of parameterized OCPs \cite{hovell_2021, art_ieeeaero24, takubo2024towards}. 
Recent efforts toward semantic conditioning include LLM-based control generation \cite{jain2026autonomous, foutter2025space-llava} and transformer models that map text to trajectories with safety enforced via SCP projection \cite{takubo2026language}. 
Yet these approaches remain tied to predefined trajectory parameterizations and waypoints, limiting their ability to encode complex mission requirements. 
Here, we instead synthesize formal mathematical constraints from text and integrate them directly into the OCP formulation.

\medskip
\noindent \textbf{LLMs for Optimization Code Generation: } 
Recent works in optimization and operations research show that LLMs are effective for extracting a structured constraint formulation from natural-language descriptions \cite{ramamonjison2023nl4opt}. 
Specifically, natural-language interfaces for translation from natural language to temporal logics \cite{wang2024chatstl, chen2023nl2tl} or automated nonconvex-to-convex transformation \cite{peng2026nc2c} are explored, with particular applications in LLM-driven safe robotic control \cite{yang2024plug, firoozi2025foundation}.
However, applications of such LLM-based constraint generation to vehicular motion planning remain limited \cite{van2024vernacopter,wang2024chatstl}, and are particularly scarce in space vehicle guidance in relation to OCPs. We address this research gap through the use of LLM–based automated constraint generation for spacecraft autonomy applications. 

\section{Methodology}

The goal of this paper is to translate natural-language mission requirements into a formal mathematical expression and then into an executable code snippet for a numerical solver, specifically a convex optimization-based algorithm. 
Formally, for a discrete-time OCP with fixed final time expressed as: 
\begin{subequations} \label{eq:ocp}
\footnotesize
\begin{alignat}{2}
    \min_{\{{x}_k, {u}_k\}_{k=1}^{N}} \quad & 
    \sum_{k=1}^{N} \mathcal{J}_k \label{eq:ocp_obj}\\
    \text{subject to} \quad 
    & {x}_{k+1} = f_k ({x}_{k} , u_k), \  
    \ && \forall k \in \mathbb{Z}_{1:N-1}, \label{eq:ocp_con_dyn} \\
    & x_1 = x_{0}, \quad  x_N = x_{f}, \label{eq:ocp_con_wyp} \\ 
    & \boldsymbol{g}(x_k, u_k) \leq 0, && \forall k \in \mathbb{Z}_{1:N-1}. \label{eq:ocp_con_ieq} \\
    & \boldsymbol{h}(x_k, u_k) = 0, && \forall k \in \mathbb{Z}_{1:N-1}, \label{eq:ocp_con_eq}
\end{alignat}
\end{subequations}
we wish to reconstruct $\boldsymbol{g}(x, u)$ and $\boldsymbol{h}(x, u)$ based on the documented mission requirements, and integrate them into an executable script. 
To this end, the proposed automated architecture is structured as a three-stage process consisting of (i) mathematical expression generation, (ii) executable code generation, and (iii) execution of a convex-based numerical OCP solver with error handling.
This intermediary step of mathematical expression generation is intended to facilitate the generation of complex requirements that require systematic decomposition with intermediate representations \cite{ma2025bridging, lyu2026canonical}.

The remainder of the section introduces the structured input design to the framework (Sec. 3.1) and the internal two-staged structure (Sec. 3.2), respectively.

\subsection{Structured input design}

Providing a structured representation of the baseline OCP formulation, together with explicit documentation of the mission requirements, is essential to enable coherent reasoning by the LLM and to support reliable code generation.
Accordingly, the proposed framework ingests three inputs:
\begin{itemize}
    \item \textbf{Base OCP code:} a Python file with the base OCP problem in code based on CVXPY library \cite{diamond2016cvxpy}, with designated placeholders indicating where additional constraints or modifications can be inserted.
    \item \textbf{Base OCP formulation:} a \LaTeX-based description of the baseline OCP that explains the general context of the problem (e.g., variable definition and base constraints).
    \item \textbf{Mission Requirement:} a text describing a feasible Scenario adjustment, trajectory requirement, or path constraint. 
\end{itemize}

\subsection{Staged constraint synthesis}

Given the structured inputs, a pre-trained LLM is queried twice with an incremental instruction towards the full constraint generation. 
Specifically, the following process is taken within the framework. 
\begin{itemize}
    \item \textbf{Mathematical Reasoning:} Instructs the LLM to augment the provided base OCP formulation in \LaTeX-format with the mathematical expression translated from the mission requirements. 
    \item\textbf{Code Generation:} Instructs LLM to read the added constraint from the generated \LaTeX-file, convert to an executable code snippet, and insert into the base OCP code at an appropriate location.
\end{itemize}

Once the executable code is generated, the numerical OCP solver is called to perform trajectory optimization, yielding either a feasible trajectory or an error report.
 
\begin{table*}[t!]
\centering
\scriptsize
\renewcommand{\arraystretch}{1.3}
\begin{tabular}{l p{0.33\linewidth}p{0.25\linewidth}p{0.25\linewidth}}
\hline
\textbf{Prompt} & \textbf{Mission Requirements} & \textbf{Mathematical Constraints} & \textbf{Example Code} \\
\hline
1 &
The target spacecraft must be observed at a distance of no more than 30m for a period of at least 15000s in the middle of the trajectory. & 
$\|(\Psi_k x_k)_{1:3}\|_2 \leq 30, \quad T_{\text{start}}\leq k\leq T_{\text{end}}$
& \texttt{[cp.norm(Psi[k][0:3]@x[k]) <= 30 for k in range(t\_start, t\_end)] }
\\
2 &
At a point in the middle of the trajectory, the relative RTN position between the two spacecraft must be [20, 0, 0] m.&
$(\Psi_k x_k)_{1:3}\ = [20,0,0] \quad k = T_{\text{mid}}$ & 
\texttt{(Psi[T\_mid]@x[T\_mid])[0:3] == cp.Constant([20,0,0])}

\\
3 &
The total Delta-V used in the first half of the trajectory must be no more than 2 mm/s.
The spacecraft must remain behind the target spacecraft in the orbit for the first 75\% of the trajectory.
& 
$\vtop{\hbox{$\begin{array}{l}
\|a\|_2 \leq 2,\quad 0\leq k\leq T_{\text{mid}} \\
(\Psi_k x_k)_{1} \leq 0,\quad 0 \leq k \leq \tfrac{3}{4} T_{\text{final}}
\end{array}$}}$
&\texttt{[cp.norm(a[k]) <=  2 for k in range(0, T\_mid)]
[Psi[T\_mid]@x[T\_mid][1] <= 0 for k in range(0, int(.75*T\_final))]
 } \\
\hline 
\end{tabular}
\caption{Mission requirements used for semantic constraint synthesis. \texttt{cp} refers to a CVXPY library \cite{diamond2016cvxpy}.}
\label{tab:problem_prompts}
\end{table*}

\section{Experiments}

A spacecraft proximity operations scenario, in which a chaser spacecraft executes a complex maneuver around a chief spacecraft, is formulated as a baseline convex OCP for the case study. 
The (relative) state representation is in Relative Orbital Elements (ROE) \cite{damico_phd_2010} scaled by the semi-major axis of the target orbit: $x_{\text{qns}} = a [\delta a, \delta \lambda, \delta e_x, \delta e_y, \delta i_x, \delta i_y]$.
Across all test cases, the initial and terminal states are fixed to $x_0 = [0, -100, 0, 0, 0, 0]$ m and $x_f = [0, 100, 0, 0, 0, 0]$ m with the chief located at a Low Earth Orbit, with a time horizon of $t_f=7.645$ hours with $N=200$ timesteps.
This corresponds to a relative reconfiguration in which the chaser initially shares the chief’s orbit with a negative along-track separation (trailing) and is required to transition to a positive along-track separation (leading).

To benchmark the performance of the proposed pipeline, three other alternative frameworks are considered: 
\begin{enumerate}[label=\Alph*)]
    \item \textbf{Single-shot, No base OCP formulation}: The pretrained model performs one-shot executable code generation using only the base OCP code and the text-based mission requirement. The base OCP formulation document is not provided.
    
    \item \textbf{Single-shot, Full inputs}: All structured inputs (base OCP formulation, base OCP code, and mission requirement) are provided to the pretrained model for one-shot executable code generation.
    \item \textbf{Two-staged, No base OCP formulation in the second stage}: The framework employs two-stage constraint synthesis via a mathematical intermediate representation. The first model (constraint translation) has access to the base OCP formulation. However, the second model (code generation) receives only (i) the base OCP code and (ii) the generated mathematical constraints.
    \item \textbf{Two-staged, Full inputs (proposed)}: The framework employs two-stage constraint synthesis via a mathematical intermediate representation. In the second stage, the code generation model is provided with (i) the base OCP code, (ii) the generated mathematical constraints, and (iii) the original base OCP formulation.
\end{enumerate}

Table~\ref{tab:problem_prompts} summarizes the documented mission requirements that are used as a part of the input, the ground truth mathematical representation, and example implementation in CVXPY \cite{diamond2016cvxpy}. 
All three requirements can be formulated as convex constraints, requiring no additional information outside what is given in the inputs in order to implement.

\subsection{Tests and results}

The experiment comprises 100 trials for each fixed input and mission requirement, evaluated across all four frameworks.
In this experiment, GPT-5.2 was adopted as a code generation model, and GPT-4.1-mini was used as a mathematical reasoning model.  

\begin{figure*}[t!]
\centering
\begin{subfigure}{0.33\textwidth}
    \centering
    \includegraphics[width=\linewidth]{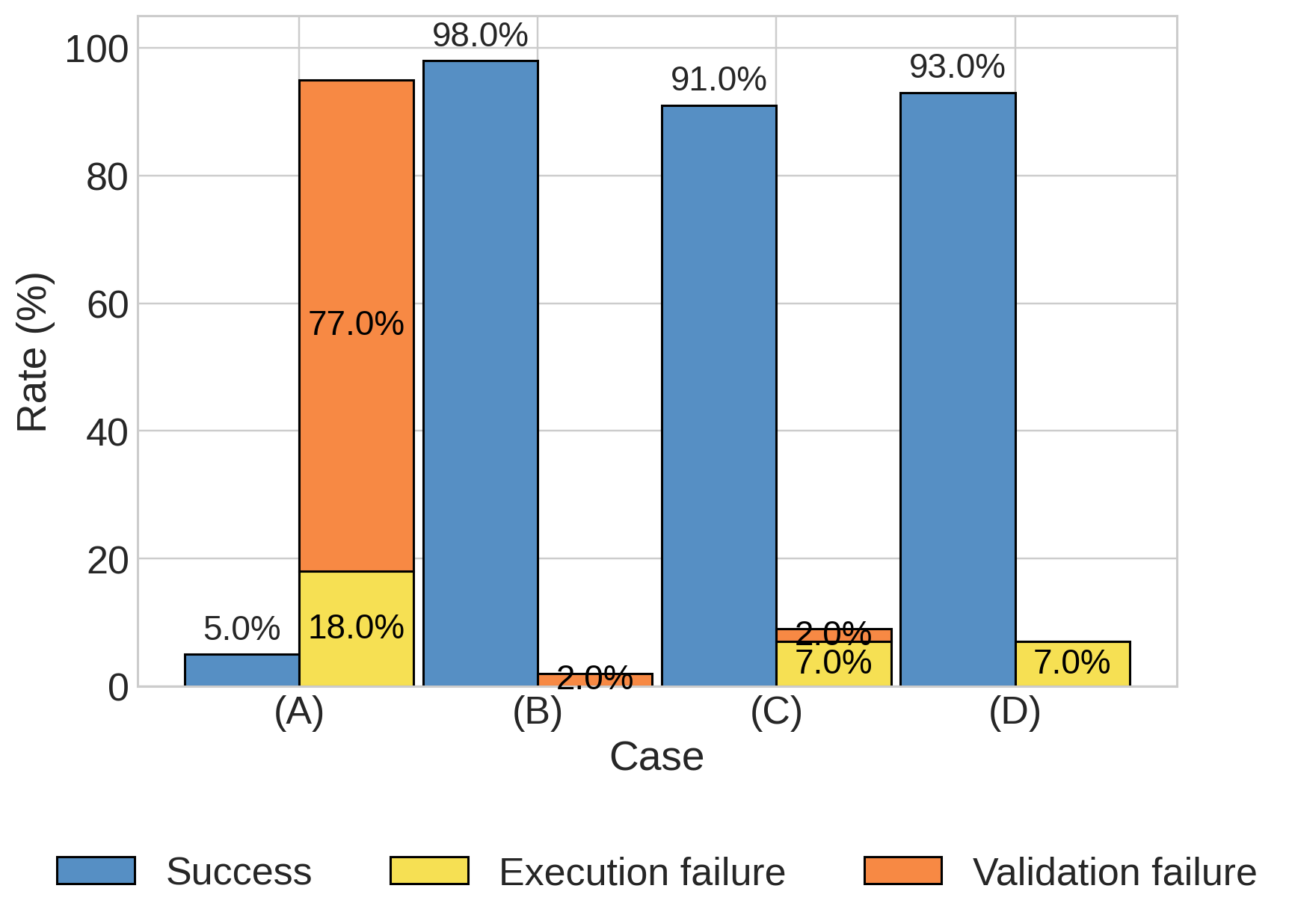}
    \caption{Prompt 1}
\end{subfigure}
\hfill
\begin{subfigure}{0.33\textwidth}
    \centering
    \includegraphics[width=\linewidth]{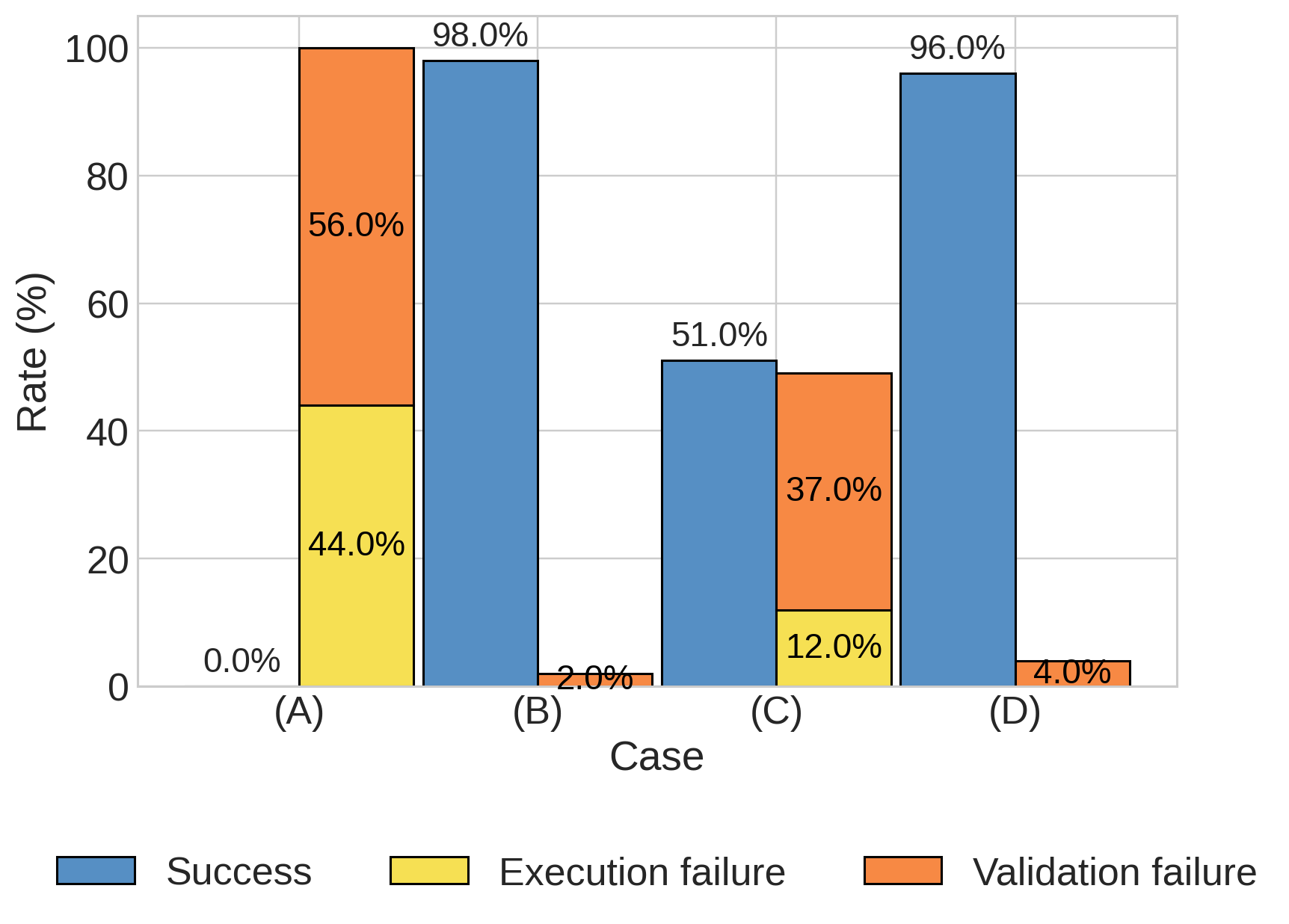}
    \caption{Prompt 2}
\end{subfigure}
\hfill
\begin{subfigure}{0.33\textwidth}
    \centering
    \includegraphics[width=\linewidth]{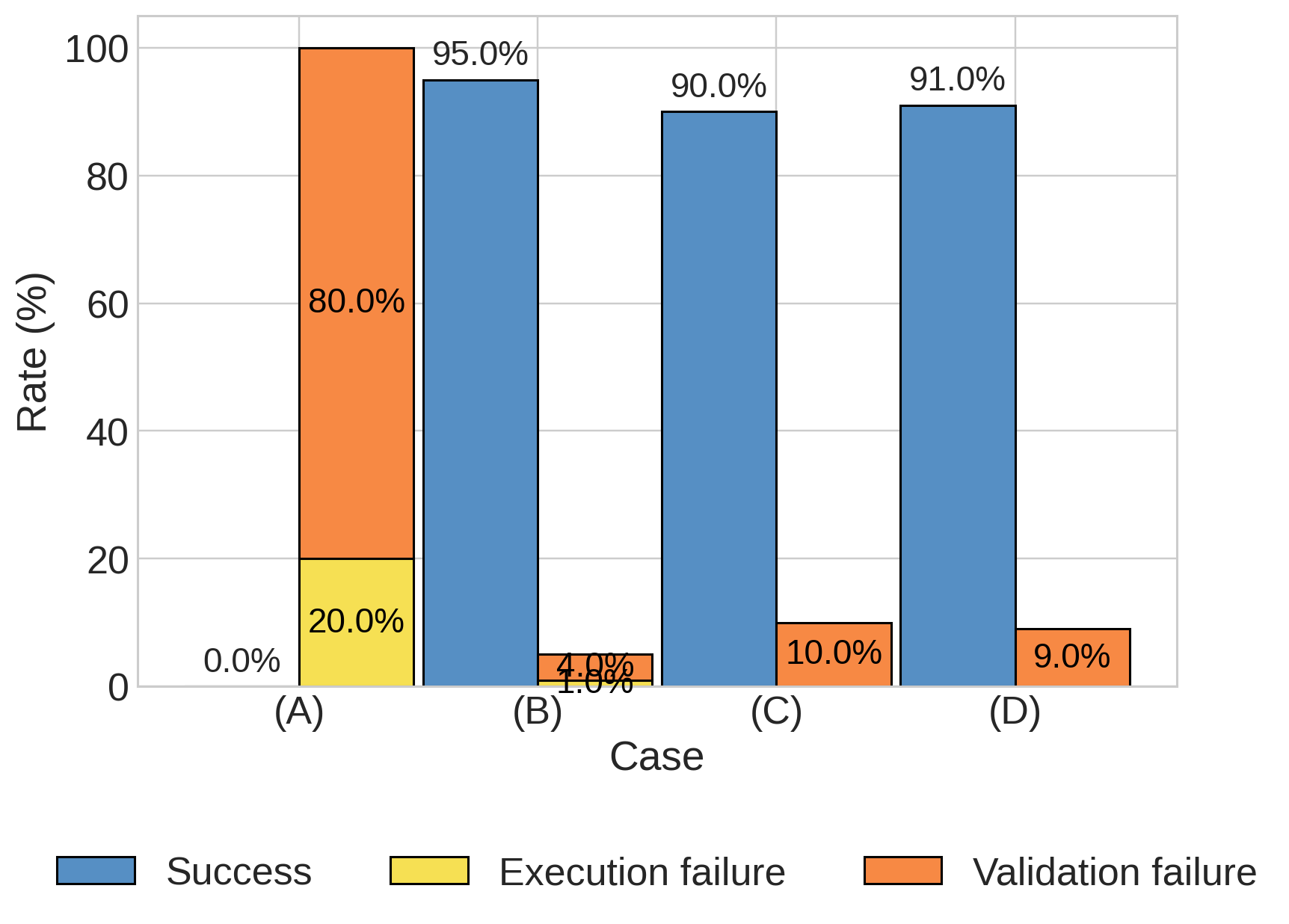}
    \caption{Prompt 3}
\end{subfigure}
\caption{Distributions of failure modes in the three prompt inputs across four frameworks: (A) Single-Shot, No OCP Formulation, (B) Single-Shot, Full Inputs, (C) Two-staged, No base OCP formulation in the second stage, and (D) Two-staged, Full inputs.}
\hfill
\label{fig:bar_graphs}
\end{figure*}

The success rates for each prompt are summarized in Table \ref{tab:success_breakdown}. 
In this paper, \textit{success} of a trial is defined as a case that generates (i) executable code and (ii) a trajectory that satisfies the target behavior specified in the mission requirements. 
Any instance in which the generated code failed to execute is classified as a failure. Furthermore, the resulting trajectory behavior must behave as expected within a margin. 
Fig.~\ref{fig:bar_graphs} summarizes the types of failures encountered with the two categories of failures documented; code execution failures and validation failures. 

\begin{table}[t!]
\centering
\footnotesize
\renewcommand{\arraystretch}{1.25}
\begin{tabular}{p{0.55\linewidth}ccc}
\hline
& \multicolumn{3}{c}{\textbf{Prompt}} \\
\cline{2-4}
& \textbf{1} & \textbf{2} & \textbf{3} \\
\hline 
(A) Single-shot, No base OCP formulation (\%) &5.0 &0.0 & 0.0\\
(B) Single-shot, Full inputs (\%) & 98.0 & 98.0 & 95.0 \\
(C) Two-staged, No base OCP formulation in the second stage (\%) & 91.0 & 51.0 & 90.0 \\
(D) Two-staged, Full inputs (\%) & 93.0 & 96.0 & 91.0 \\
\hline 
\end{tabular}
\caption{Success statistics across three prompts (100 trials each).}
\label{tab:success_breakdown}
\end{table}

From these results, the importance of a structured input that provides the problem context is clear. 
In cases where there is no initial base OCP formulation, the model almost always fails to get the desired trajectory behavior, due to a mix of code-based failures and validation failures (cf. Fig.~\ref{fig:bar_graphs}). 
These results underscore the importance of providing a structured input that provides sufficient problem context. 
In contrast, for scenarios where the OCP formulation is given, the success rate of the single-shot code generation pipeline is markedly improved, with success rates exceeding 95\%  across all tested prompts. 

When comparing the two-stage frameworks and the effect of including the full inputs at each stage, the success rate was approximately maintained for Prompt 1 and 3 while Prompt 3 suffered a significant drop. This suggests that the OCP formulation provides essential contextual information, and passing it to the second stage helps prevent failures caused by subtle but crucial missing details that may be omitted from the generated math formulation.

While all frameworks that incorporate a base OCP formulation at each stage achieve success rates exceeding 90\% across the evaluated prompts, the single-shot approach with full inputs demonstrates the strongest performance under the selected mission requirements. This result may be attributed to the relatively simple nature of the tested requirements, all of which have direct convex implementations. If this interpretation is correct, the proposed Two-Stage method—despite maintaining a high success rate—offers limited advantage over a more direct formulation strategy in these cases. However, as requirement complexity increases and the prompt-to-code pathway requires additional intermediate reasoning steps, we expect the performance of the direct single-shot method to degrade. Given that the long-term objective of this framework is to address such complex requirements and generate their corresponding mathematical formulations, the Two-Stage architecture remains the more promising approach for scalable and flexible automated constraint generation.

\section{Conclusion and Future Work}

This paper proposes an automated framework that translates semantic mission requirements into well-defined mathematical constraints and executable code for trajectory optimization problems.
The results indicate that consistent generation of executable code from semantic mission requirements is achievable with appropriate input structuring and prompt design.
In the evaluated prompts, an intermediate mathematical formulation step was not strictly necessary; however, incorporating this step did not significantly degrade performance.
Future work will extend the framework to (i) an automatic convexification with formal verification to accommodate nonconvex constraints, (ii) ablation studies across different pretrained LLMs, and (iii) supervised fine-tuning of the models and associated scalable dataset curation tailored for trajectory optimization.

{
    \small
    \bibliographystyle{ieeenat_fullname}
    \bibliography{main}
}

\clearpage
\setcounter{page}{1}
\maketitlesupplementary

\section{Generation Prompts}
\begin{promptbox} {\LaTeX  \quad Generation Prompt Example}
You are an expert in optimal control problem formulation and convex optimization using cvxpy.

You will take in 2 files:

1. Problem Spec: LATEX file that lays out a given optimal control problem
2. Problem Prompt: A description of a change to the given control problem in that will result in a new constraint or set of constraints for the problem 

Your job is to output a new Latex document with a section for each new constraint that includes:
-The minimal mathematical formulation of the constraint in terms of the variables and parameters defined in the original problem spec.
- A single line mathematical equation for the constraint in terms of the variables and parameters defined in the original problem spec and any values that you have also defined in the mathematical formulation of the constraint.

Rules for the constraint addition section:
- Do not use any explanatory text, only include the mathematical formulation and the single line mathematical equation.
-  Ensure the Dimensions of all variables is explicitly defined in the mathematical formulation of the constraint. 
- This description will be used to write code so ensure that the mathematical formulation is clear and unambiguous.


Output nothing other than the LATEX document with just the constraints. Do not include any explanations or commentary.
\end{promptbox}

\begin{promptbox} {Code Generation Prompt Example}
You are an expert in optimal control problem formulation and convex optimization using cvxpy.
You are editing an existing Python file.
Based on the provided constraint addition latex file and the latex spec file you will add the constraints in the constraint addition latex file.

Rules:
- You will keep all existing code.
- ONLY replace code between:
  # BEGIN AUTO-GENERATED CONSTRAINTS
  and
  # END AUTO-GENERATED CONSTRAINTS
- Do NOT change anything else.
- Return the FULL updated file.
- No markdown. No explanation.



Generate constraint code that:
- Only uses the information explicitly provided in the provided specification and does not require further user input to define the constraint. 
- Ensures that the constraint is satisfied when discretized.
- Is mindful of dimensionality and ensures all variables are defined with explicit dimensions in the code.
- Does not overly complicate the constraint and uses the simplest form that satisfies the requirement of the constraint.
- Uses cvxpy syntax
- Is DCP compliant

\end{promptbox}

\section{Example Output}

\subsection{Generated \LaTeX{} Output for Prompt 1}

\small 
\begin{framed}
\subsection*{Observation Distance Constraint}

Let $d_{\max} = 30$ m be the maximum allowable relative position norm for observation.

Let $T_{\text{obs}} = 15000$ s be the required observation duration.

Let $\Delta t$ be the time step duration in seconds.

Define the number of discrete time steps for observation as
\[
M = \left\lceil \frac{T_{\text{obs}}}{\Delta t} \right\rceil, \quad M \in \mathbb{N}.
\]

Let $t_{\text{start}} = \left\lfloor \frac{N - M}{2} \right\rfloor$ be the starting time index of the observation window.

For all $t = t_{\text{start}}, t_{\text{start}}+1, \ldots, t_{\text{start}} + M - 1$, the relative position vector in RTN frame is
\[
\mathbf{p}_t = \left( \mathbf{x}^{\text{RTN}}_t \right)_{1:3} = \left( \Psi_t \mathbf{s}_t \right)_{1:3} \in \mathbb{R}^3.
\]

The constraint is
\[
\|\mathbf{p}_t\|_2 \leq d_{\max}, \quad \forall t = t_{\text{start}}, \ldots, t_{\text{start}} + M - 1.
\]

\bigskip

\noindent\textbf{Single line constraint:}
$$\|\left(\Psi_t \mathbf{s}_t\right)_{1:3}\|_2 \leq 30 $$
$$ t = \left\lfloor \frac{N - \lceil T_{obs}/\Delta t \rceil}{2} \right\rfloor, \ldots, \left\lfloor \frac{N - \lceil T_{obs}/\Delta t \rceil}{2} \right\rfloor $$

$$+ \lceil T_{obs}/\Delta t \rceil - 1 $$

\end{framed}
\normalsize
\newpage

\subsection{Generated Code Output Snippet For Prompt 1}

\begin{lstlisting}[language=python, breaklines=true, basicstyle=\ttfamily\scriptsize, frame=single] 
#BEGIN AUTO-GENERATED CONSTRAINTS####################

# Observation Distance Constraint:
# For an observation window of duration T_obs=15000s, enforce ||p_t||_2 <= d_max=30m
# where p_t = (Psi_t s_t)_{1:3}.
d_max = 30.0
T_obs = 15000.0

# Conservative discretization: ensure at least one step if dt is large; cap to horizon length.
M = int(np.ceil(T_obs / float(self.dt))) if self.dt is not None else self.n_time
M = max(1, M)
M = min(M, self.n_time)

t_start = int(np.floor((self.n_time - M) / 2.0))
t_start = max(0, min(t_start, self.n_time - M))

for t in range(t_start, t_start + M):
    # psi[t] is expected shape (6,6); s[t] is (6,)
    p_t = (psi[t] @ s[t])[0:3]
    con += [cp.norm(p_t, 2) <= d_max]    
#END AUTO-GENERATED CONSTRAINTS####################
\end{lstlisting}

\end{document}